\documentstyle[12pt]{article}
 \normalsize

\begin{document}
\title{{\bf {On asymptotic properties of the generalized Dirichlet $L$-functions}}
\footnotetext{This work is supported by Natural Science Basic
Research Program of Shanxi (2016JQ1013), Fundamental Research
Funds for the Central Universities
of Northwestern Polytechnical University(3102016ZY030) and Industrial and Educational Cooperation and Collaborative Education Project of the Ministry of Education(201801332030).\\
1.marong@nwpu.edu.cn \qquad 2. 1332524636@qq.com\qquad 3. zzboyzyl@163.com}
\author{Rong
Ma$^{1}$   \\{\small{{School of Science, Northwestern
Polytechnical University}} }\\{\small{Xi'an, Shaanxi, 710072,
People's Republic of China}}\\{\small{{School of Mathematics and Statistics, University of Glasgow}} }\\{\small{Glasgow, G128QQ,
United Kingdom}}\\ Yana Niu$^{2}$  \\{\small{{School of Science, Northwestern
Polytechnical University}} }\\{\small{Xi'an, Shaanxi, 710072,
People's Republic of China}} \\ Yulong
Zhang$^{3}$\\{\small{{School of Software Engineering, Xi'an
Jiaotong University}}}
\\{\small{Xi'an, Shaanxi, 710049, People's Republic of China}}\\}
\date{}}
\maketitle
\vspace{0.0cm}

\begin{center}
\large{{\bf Abstract}}
\end{center}

\small{Let $q\ge3$ be an integer, $\chi$ denote a
Dirichlet character modulo $q$, for any real number $a\ge 0$, we
define the generalized Dirichlet $L$-functions
$$
L(s,\chi,a)=\sum_{n=1}^{\infty}\frac{\chi(n)}{(n+a)^s},
$$
where $s=\sigma+it$ with $\sigma>1$ and $t$ both real. It can be
extended to all $s$ by analytic continuation. In this paper, we
study the mean value properties of the generalized Dirichlet
$L$-functions, and obtain several
sharp asymptotic
formulae by using analytic method.
}

\vspace{1cm}
\begin{center}
\large{{\bf1. Introduction}}
\end{center}

Let $q\ge3$ be an integer, $\chi$ denote a Dirichlet character modulo
$q$, Dirichlet $L$-functions
$L(s,\chi)$ defined by
$$
L(s,\chi)=\sum_{n=1}^{\infty}\frac{\chi(n)}{n^s},
$$
where $s=\sigma+it$ with $\sigma>1$
and $t$ both real. It is very important in analytic number theory, and many studies have be done in all directions of Dirichlet $L$-functions. One of the most significant aspects about Dirichlet $L$-functions are the mean value properties. D.R. Heath-Brown, W. Zhang, R. Balasubramanian(see Ref. [1-3]) studied the square
mean value properties on Dirichlet $L$-functions on the line
$\sigma=\frac{1}{2}$. For example, R. Balasubramanian (see Ref. [3]) got
the asymptotic formula
$$
\sum_{\chi\bmod{q}}\left|L\left(\frac{1}{2}+it,\chi\right)\right|^2=\frac{\phi^2(q)}{q}\log(qt)+O(q(\log\log
q)^2)+O(te^{10\sqrt{\log
q}})+O(q^{\frac{1}{2}}t^{\frac{2}{3}}e^{10\sqrt{\log q}}),
$$
which is satisfied for $t\ge3$ and for all $q$.

W. Zhang, Y. Yi(see Ref. [4] and [5]) got
different kinds of the mean value of Dirichlet $L$-functions with
weight or not. For instance, Y. Yi and W. Zhang (see
Ref.\,[5]) gave the asymptotic formula of Dirichlet $L$-functions
with the weight of $\tau(\chi)$
\begin{eqnarray}
&&\sum_{\chi\ne\chi_0}|\tau(\chi)|^m|L(1,\chi)|^{2k}\nonumber\\
&=&N^{\frac{m}{2}-1}\phi^{2}(N)\zeta^{2k-1}(2)\prod_{p\mid
q}(1-\frac{1}{p^2})^{2k-1}\prod_{p\dag
q}(1-\frac{1-C_{2k-2}^{k-1}}{p^2})\prod_{p\mid
M}(p^{\frac{m}{2}+1}-2p^{\frac{m}{2}}+1)\nonumber\\
&&+O(q^{\frac{m}{2}+\epsilon}),\nonumber
\end{eqnarray}
where $q$ is an integer $\ge3$ and
$q=MN,(M,N)=1,M=\prod\limits_{p|q \atop p^2\dag q}p$.

The first and third authors (see Ref. [6]) also studied the mean value of Dirichlet $L$-functions with trigonometric sums, and gave the asymptotic formula as follow \\
\begin{eqnarray}
&&\sum_{\chi\ne\chi_0}\left|\sum_{a=1}^{p-1}\chi(a)e\left(\frac{f(a)}{p}\right)\right|^2
|L(1,\chi)|^{2m}\nonumber\\
&=&p^2\,\zeta^{2m-1}(2)
\prod_{p_0}\lgroup
1-\frac{1-C_{2m-2}^{m-1}}{p_0^2}\rgroup+O(p^{2-\frac{1}{k}+\epsilon}
),\nonumber
\end{eqnarray}
where p is prime, $\epsilon$ is any small positive real number, $f(x)=\sum_{i=0}^ka_ix^i$ is a polynomial such that deg$(f(x))=k$
and $p\dag (a_0,a_1,\dots,a_k)$, $m$
and $k$ are any positive integers, $\prod_{p_0}$ denotes the product over all primes different
from $p$, $C_{m}^{n}=\frac{m!}{n!(m-n)!}$, and the $O$ constant
depends only on $k$ and $\epsilon$. Obiviously, let $f(a)=a$, we have
the mean value of Dirichlet
$L$-functions with the weight of $\tau(\chi)$.

Now let $a\ge0$ be an integer, generalized Dirichlet $L$-functions
$L(s,\chi,a)$ defined by
$$
L(s,\chi,a)=\sum_{n=1}^{\infty}\frac{\chi(n)}{(n+a)^s},
$$
where $s=\sigma+it$ with $\sigma>1$
and $t$ both real.

About the generalized Dirichlet series, B.
C. Berndt (see Ref. [7]-[9]) studied many identical properties
satisfying restrictive conditions. It is well known that for
$\chi$ a nonprincipal, primitive character modulo $q$, for
$\sigma>\frac{1}{2}-m$ with $m$ a positive integer, Prof. Berndt
(see Ref. [9]) derived
$$
L(s,\chi,a)=\frac{a^{-s}}{\Gamma(s)}\left(\sum_{j=0}^{m-1}\frac{(-1)^{j}\Gamma(s+j)L(-j,\chi)}{j!\,a^{j}}+G(s)\right),
$$
where $G(s)$ is an analytic function. When $n$ is a nonpositive
integer, we can easily calculate $L(n,\chi,a)$, in particular,
$L(0,\chi,a)=L(0,\chi)$.

The first and third authors (see Ref. [10]) also got the
following asymptotic formula about the generalized Dirichlet
$L$-functions
\begin{eqnarray}
\sum_{\chi\ne\chi_0}\left|L(1,\chi,a)\right|^2
=\phi(q)\sum_{d|q}\frac{\mu(d)}{d^2}\zeta\left(2,\frac{a}{d}\right)
-\frac{4\phi(q)}{a}\sum_{d|q}\frac{\mu(d)}{d}\sum_{k=1}^{[\frac{a}{d}]}\frac{1}{k}
+O\left(\frac{\phi(q)\log q}{\sqrt{q}}\right),
\end{eqnarray}
where $\zeta(s,\alpha)(s=\sigma+it,\,\alpha>0)$ is the Hurwitz
zeta function
 defined for
$\sigma>1$ by the series
$$
\zeta(s,\alpha)=\sum_{n=0}^{\infty}\frac{1}{(n+\alpha)^{s}},
$$
and $\phi$ is the Euler function, $\mu$ is the M\"{o}bius
function, the $O$ constant only depends on $a$.

On the other hand, trigonometric sums are also
the most important research topic in analytic number theory. Let
$p$ be a prime, $f(x)=a_0+a_1x+\dots +a_kx^k$ is a $k$-degree
polynomial with integral coefficients such that
$(p,a_0,a_1,\dots,a_k)=1$, trigonometric sums are defined by
$$
\sum_{a=1}^{p}\chi(a)e\left(\frac{f(a)}{p}\right),
$$
where $\chi$ denotes a Dirichlet Character modulo $p$ and
$p\dag (a_0,a_1,\dots,a_k)$. When $\chi=\chi_0$, we can see the
trigonometric sums
enjoy many good properties (see Ref. [11-15]).

If communicated
with the generalized trigonometric sums or with Dirichlet character, whether the generalized
Dirichlet $L$-functions still show good properties? The authors are
very interested in the problems. But there are few references to
be referred to about these problems. In this paper, we will study
the mean value properties of the generalized Dirichlet $L$-functions
with Dirichlet character and the generalized trigonometric sums,
$$
\sum_{\chi\ne\chi_0}\chi(k)
|L(1,\chi,a)|^{2},
$$
$$
\sum_{\chi\ne\chi_0}\left|\sum_{x=1}^{p-1}\chi(x)e\left(\frac{f(x)}{p}\right)\right|^2
|L(1,\chi,a)|^{2},
$$
where $p\ge3$ is odd prime,$\chi$ is a Dirichlet character modulo
$p$, $\chi_0$ is the non-principal character modulo $p$, $a\ge 0$
is any real number with $(a,p)=1$, $e(y)=e^{2\pi iy}$. It could
tell us some relationship between the Dirichlet
character and the generalized trigonometric sums. More
precisely, we prove the following theorems:

\noindent \textbf{Theorem 1. }Let $k, q$ be two integers with
$ q\ge 3, k\not=1, (k,q)=1$ and $\chi$ denote a Dirichlet character modulo $q$.
Then for any positive real number $a\ge1$, we have the asymptotic
formula
\begin{eqnarray}
\sum_{\chi\ne\chi_0}\chi(k) |L(1,\chi,a)|^{2}
=\frac{\phi(q)}{a(k-1)}\sum_{d|q}\frac{\mu(d)}{d}\sum_{l=[\frac{a}{kd}]+1}^{[\frac{a}{d}]}\frac{1}{l}
+O\left(\frac{\phi(q)\log q}{\sqrt{q}}\right),\nonumber
\end{eqnarray}
where $\phi$ is the Euler function, $\mu$ is the
M\"{o}bius function, and the $O$ constant depends only on $a,k$.

\noindent \textbf{Theorem 2. }Let $p\ge 3$ be an odd prime and
$\chi$ denote a Dirichlet character modulo $p$,
$f(x)=a_0+a_1x+\dots +a_kx^k$ is a $k$-degree polynomial with
integral coefficients such that $(p,a_0,a_1,\dots,a_k)=1$,
$e(y)=e^{2\pi iy}$. Then for any positive real number $a\ge1$ with
$(a,p)=1$, we have the asymptotic formula
\begin{eqnarray}
&&\sum_{\chi\ne\chi_0}\left|\sum_{x=1}^{p-1}\chi(x)e\left(\frac{f(x)}{p}\right)\right|^2
|L(1,\chi,a)|^{2}\nonumber\\
&=&p^2\sum_{d|q}\frac{\mu(d)}{d^2}\zeta\left(2,\frac{a}{d}\right)
-\frac{4p^2}{a}\sum_{d|q}\frac{\mu(d)}{d}\sum_{l=1}^{[\frac{a}{d}]}\frac{1}{l}
+O(p^{2-\frac{1}{k}+\epsilon}),\nonumber
\end{eqnarray}
where $\mu$ is the M\"{o}bius function and
$\zeta(s,\alpha)(s=\sigma+it,\,\alpha>0)$ is the Hurwitz zeta
function. The $O$ constant is depending on $k$,$a$ and $\epsilon$.

\textbf{Note 1.} For the general case of $2l$-th ($l\ge2$) power mean value of the
generalized Dirichlet $L$-functions and the Dirichlet charater
$$
\sum_{\chi\ne\chi_0}\chi(k) |L(1,\chi,a)|^{2l},
$$
it is still an open problem.

\textbf{Note 2.} For the general case of $2m$-th ($m\ge2$) power of the generalized
trigonometric sums and $2l$-th ($l\ge2$)power mean value of the
generalized Dirichlet $L$-functions
$$
\sum_{\chi\ne\chi_0}\left|\sum_{x=1}^{p-1}\chi(x)e\left(\frac{f(x)}{p}\right)\right|^{2m}
|L(1,\chi,a)|^{2l},
$$
it is still an open problem.

\begin{center}
\large{{\bf 2. Some lemmas}}
\end{center}

To complete the proofs of both of the Theorems,\,we need the following several
lemmas.\,First, we make an identity of
the Dirichlet $L$-functions and the generalized form.

\noindent \textbf{Lemma 1. }Let $q\ge3$ be an integer, and $\chi$
denote a nonprincipal Dirichlet character modulo $q$. Let $L(s,\chi)$ denote the
Dirichlet $L$-functions corresponding to $\chi$, and $L(s,\chi,a)$
denote the generalized Dirichlet $L$-functions. Then for any real
number $a\ge0$, we have
$$
L(1,\chi,a)=L(1,\chi)-a\sum_{n=1}^{\infty}\frac{\chi(n)}{n(n+a)}.
$$
\textbf{Proof. }See Lemma 1, and let $m=1$ (Ref. [10]).

\noindent \textbf{Lemma 2. }Let $f(x)$  be a polynomial with
integer coefficients as $ f(x)=a_0+a_1x+\dots+a_kx^k,$\,and $\chi$
be a Dirichlet character modulo $p$. Then
we have\\
$$
\left|\sum_{x=1}^{p-1}\chi(x)e\left(\frac{f(x)}{p}\right)\right|^2=p-1+
\sum_{x=2}^{p-1}\chi(x)\sum_{y=1}^{p-1}e\left(\frac{g(y,x)}{p}\right),
$$
where $g(y,x)=f(xy)-f(y)=\sum_{i=0}^{k}a_i(x^i-1)y^i$.\\
\textbf{Proof. }Note that for $1\le y\le p-1$ ($p$ is prime), we have $(y,p)=1$.
According to the properties of characters, we have\\
\begin{eqnarray}
\left|\sum_{x=1}^{p-1}\chi(x)e\left(\frac{f(x)}{p}\right)\right|^2
&=&\sum_{x,\,y=1}^{p-1}\chi(x)\bar{\chi}(y)e\left(\frac{f(x)-f(y)}{p}\right)\nonumber\\
&=&\sum_{x=1}^{p-1}\sum_{y=1}^{p-1}\chi(xy)\bar{\chi}(y)
e\left(\frac{f(xy)-f(y)}{p}\right).\nonumber
\end{eqnarray}
Let \\
$$g(y,x)=f(xy)-f(y)=\sum_{i=0}^{k}a_i(x^i-1)y^i,$$
we get\\
\begin{eqnarray}
\left|\sum_{x=1}^{p-1}\chi(x)e\left(\frac{f(x)}{p}\right)\right|^2
&=&\sum_{x=1}^{p-1}\sum_{y=1}^{p-1}\chi(x)e\left(\frac{g(y,x)}{p}\right)\nonumber\\
&=&\sum_{y=1}^{p-1}\chi(1)e\left(\frac{g(y,1)}{p}\right)+
\sum_{x=2}^{p-1}\sum_{y=1}^{p-1}\chi(x)e\left(\frac{g(y,x)}{p}\right)\nonumber\\
&=&p-1+\sum_{x=2}^{p-1}\chi(x)\sum_{y=1}^{p-1}e\left(\frac{g(y,x)}{p}\right).\nonumber
\end{eqnarray}
This proves Lemma 2. \\

\noindent \textbf{Lemma 3. }Let $f(x)$ satisfy the conditions of
Lemma 2. And let
$g(z)=g(z,x)=f(xz)-f(z)=\sum_{i=0}^{k}a_i(x^i-1)z^i$,\,then we
have the following
estimate\\
$$
\left|\sum_{y=1}^{p-1}e\left(\frac{g(y,x)}{p}\right)\right|\,\,\,\,\left\{
\begin{array}{ll}
\ll p^{1-\frac{1}{k}},& p\dag(b_0,b_1,\dots,b_k)  \\
=p-1,& p|(b_0,b_1,\dots,b_k)
\end{array}
\right.
$$
where $b_i=a_i(x^i-1),i=0,1,\dots,k$ and $k$ is the degree of the polynomial $f(x)$.\\
\textbf{Proof. }The result is apparent if
$p\,|(b_0,b_1,\dots,b_k)$.\,If $p\dag(b_0,b_1,\dots,b_k)$,\,
according to the definition of $g(x,a)$,\,we have (see Ref. [12]),
$$
\left|\sum_{b=1}^{p-1}e\left(\frac{g(b,a)}{p}\right)\right|\ll
p^{1-\frac{1}{k}}.
$$
This proves Lemma 3.\\

\noindent \textbf{Lemma 4. }Let $q\ge 3$ be an integer and $\chi$
be the Dirichlet character modulo $q$.\,Then for any positive integer $a\ge2$ with $(a,q)=1$, we have\\
$$
\sum_{\chi\ne \chi_{0}}\chi(a)\arrowvert L(1,\chi)
\arrowvert^{2}=\frac{\phi(q)}{a}\zeta(2)\prod_{p|q}\left(1-\frac{1}{p^2}\right)+O(\log^2
q).
$$
\textbf{Proof. }For convenience, we put \\
$$
A(\chi,y)=\sum_{\frac{q}{a}\le n\le
y}\chi(n),\,\,\,B(\chi,y)=\sum_{q\le n\le y}\chi(n).
$$
Then for $s>1$,\,the series $L(s,\chi)$ is absolutely convergent,
so applying Abel's identity we have
\begin{eqnarray}
L(s,\chi)&=&\sum_{n=1}^{\infty}\frac{\chi(n)}{n^s}\nonumber\\
&=&\sum_{n=1
}^{\frac{q}{a}}\frac{\chi(n)}{n^s}+s\int_{\frac{q}{a}}^{+\infty}\frac{A(\chi,y)}{y^{s+1}}dy\nonumber\\
&=&\sum_{n=1}^{q}\frac{\chi(n)}{n^s}+s\int_{q}^{+\infty}\frac{B(\chi,y)}{y^{s+1}}dy.\nonumber
\end{eqnarray}

It is clear that the above formula also holds for $s=1$ and $\chi
\ne \chi_{0}$. Hence according to the definition of the Dirichlet
$L$-function, for any positive integer $a\not=1$ and $(a,q)=1$, we
have
\begin{eqnarray}
& & \sum_{\chi\ne \chi_{0}}\chi(a)\arrowvert L(1,\chi)
\arrowvert^{2}\nonumber\\
&=&\sum_{\chi\ne \chi_{0}}\chi(a)\left|
\sum_{n=1}^{\infty}\frac{\chi(n)}{n}
\right|^{2}\nonumber\\
&=&\sum_{\chi\ne \chi_{0}}\chi(a)\left(
\sum_{n=1}^{\infty}\frac{\chi(n)}{n}\right)\left(\sum_{l=1}^{\infty}\frac{\bar{\chi}(l)}{l}
\right)\nonumber\\
&=&\sum_{\chi\ne \chi_{0}}\chi(a)\left( \sum_{1\le n\le
\frac{q}{a}}\frac{\chi(n)}{n}+\int_{\frac{q}{a}}^{+\infty}\frac{A(\chi,y)}{y^{2}}dy
\right)\left(
\sum_{l=1}^{q}\frac{\bar{\chi}(l)}{l}+\int_{q}^{+\infty}\frac{B(\bar{\chi},y)}{y^{2}}dy
\right)\nonumber\\
&=&\sum_{\chi\ne \chi_{0}}\chi(a)\left( \sum_{n=1
}^{\frac{q}{a}}\frac{\chi(n)}{n} \right)\left(
\sum_{l=1}^{q}\frac{\bar{\chi}(l)}{l}
\right)\nonumber\\
& &+\sum_{\chi\ne \chi_{0}}\chi(a)\left( \sum_{n=1
}^{\frac{q}{a}}\frac{\chi(n)}{n} \right)\left(
\int_{q}^{+\infty}\frac{B(\bar{\chi},y)}{y^{2}}dy
\right)\nonumber\\
& &+\sum_{\chi\ne \chi_{0}}\chi(a)\left(
\sum_{l=1}^{q}\frac{\bar{\chi}(l)}{l}\right)\left(
\int_{\frac{q}{a}}^{+\infty}\frac{A(\chi,y)}{y^{2}}dy
\right)\nonumber\\
& &+\sum_{\chi\ne
\chi_{0}}\chi(a)\left(\int_{\frac{q}{a}}^{+\infty}
\frac{A(\chi,y)}{y^{2}}dy\right)\left(\int_{q}^{+\infty}\frac{B(\bar{\chi},y)}{y^{2}}dy\right)\nonumber\\
&\equiv& A_{1}+A_{2}+A_{3}+A_{4}.\nonumber
\end{eqnarray}
Now we will estimate each of the above term .\\

From the orthogonality relation for character sums
    modulo $q$,\,we know that for $(q,n)=1$, \,we have the identity\\
    $$
\sum_{\chi \bmod q}\chi(n)\bar{\chi}(l)=\left\{
\begin{array}{ll}
\phi(q),&\mbox{if $n\equiv l$ mod $q$ };\\
0,&      \mbox{otherwise}.
\end{array}
\right.
    $$\\
Then we can easily get \\
\begin{eqnarray}
A_{1}&=&\sum_{\chi\ne \chi_{0}}\chi(a)\left( \sum_{n=1
}^{\frac{q}{a}}\frac{\chi(n)}{n} \right)\left(
\sum_{l=1}^{q}\frac{\bar{\chi}(l)}{l}
\right)\nonumber\\
&=&\sum_{\chi \bmod q}\,\,\sum_{n=1}
^{\frac{q}{a}}\sum_{l=1}^{q}\frac{\chi(an)\bar{\chi}(l)}{nl}-\sum_{n=1}
^{\frac{q}{a}}\sum_{l=1}^{q}\frac{1}{nl}\nonumber\\
&=&\phi(q)\displaystyle\mathop{\displaystyle\mathop{\sum{'}}_{n=1}^{\frac{q}{a}}
\displaystyle\mathop{\sum{'}}_{l=1}^{q}}_{an\equiv l(\bmod q)}\frac{1}{nl}+O(\log^2 q)\nonumber\\
&=&\phi(q)\displaystyle\mathop{\sum{'}}_{n=1}^{\frac{q}{a}}\frac{1}{an^2}+O(\log^2 q)\nonumber\\
&=&\phi(q)\sum_{n=1}^{\frac{q}{a}}\frac{1}{an^2}\sum_{d|(n,q)}\mu(d)+O(\log^2 q)\nonumber\\
&=&\phi(q)\sum_{d|q}\mu(d)\sum_{n=1\atop
d|n}^{\frac{q}{a}}\frac{1}{an^2}+O(\log^2 q)\nonumber\\
&=&\frac{\phi(q)}{a}\sum_{d|q}\frac{\mu(d)}{d^2}\sum_{n=1}^{\frac{q}{ad}}\frac{1}{n^2}+O(\log^2
q)\nonumber\\
&=&\frac{\phi(q)}{a}\sum_{d|q}\frac{\mu(d)}{d^2}\sum_{n=1}^{\infty}\frac{1}{n^2}
+O\left(\frac{\phi(q)}{a}\sum_{d|q}\frac{\mu(d)}{d^2}\sum_{n=q/ad}^{\infty}\frac{1}{n^2}\right)+O(\log^2
q)\nonumber\\
&=&\frac{\phi(q)}{a}\zeta(2)\prod_{p|q}\left(1-\frac{1}{p^2}\right)+O(\log^2
q),
\end{eqnarray}
where $\sum{'}_{n}$ indicates that the sum is over those $n$ relatively prime to $q$.

According to Cauchy inequality and Polya-Vinogradiv inequality about character sums we can easily get
\begin{eqnarray}
A_{2}&=&\sum_{\chi\ne \chi_{0}}\chi(a)\left( \sum_{n=1
}^{\frac{q}{a}}\frac{\chi(n)}{n} \right)\left(
\int_{q}^{+\infty}\frac{B(\bar{\chi},y)}{y^{2}}dy
\right)\nonumber\\
&=&\sum_{\chi\ne \chi_{0}}\chi(a)\left( \sum_{n=1
}^{\frac{q}{a}}\frac{\chi(n)}{n} \right)\left(
\int_{q}^{q^{\frac{3}{2}}}\frac{\sum_{q\le n\le
y}\bar{\chi}(n)}{y^{2}}dy
\right)\nonumber\\
& & +\sum_{\chi\ne \chi_{0}}\chi(a)\left( \sum_{n=1
}^{\frac{q}{a}}\frac{\chi(n)}{n} \right)\left(
\int_{q^{\frac{3}{2}}}^{+\infty}\frac{\sum_{q\le n\le
y}\bar{\chi}(n)}{y^{2}}dy
\right)\nonumber\\
&\le&\int_{q}^{q^{\frac{3}{2}}}\frac{1}{y^2}\left|\sum_{n=1}^{\frac{q}{a}}
\sum_{l=q}^{y}\frac{1}{n}\sum_{\chi\ne
\chi_{0}}\chi(an)\bar{\chi}(l)\right|
dy+q^{\epsilon}\int_{q^{\frac{3}{2}}}^{+\infty}\frac{1}{y^2}
\sum_{\chi\ne \chi_{0}}\arrowvert B(\bar{\chi},y)\arrowvert
dy\nonumber\\
&\ll&\int_{q}^{q^{\frac{3}{2}}}\frac{\phi(q)}{y^2}\left|
\displaystyle\mathop{\displaystyle\mathop{\sum{'}}_{n=1}^{\frac{q}{a}}
\displaystyle\mathop{\sum{'}}_{l=q}^{y}}_{\atop an\equiv l(\bmod
q)}\frac{1}{n}\right|
dy+q^{\epsilon}\int_{q^{\frac{3}{2}}}^{+\infty}\frac{1}{y^2}\left(\sum_{\chi\ne
\chi_{0}}1^2\right)^{\frac{1}{2}}\left(\sum_{\chi\ne
\chi_{0}}\arrowvert
B(\bar{\chi},y)\arrowvert^2\right)^{\frac{1}{2}}dy\nonumber\\
&\le&\phi(q)\int_{q}^{q^{\frac{3}{2}}}\frac{1}{y^2}\frac{y}{q}\displaystyle\mathop{\sum{'}}_{n=1}^{\frac{q}{a}}\frac{1}{n}dy
+\phi(q)q^{\frac{1}{2}+\epsilon}\int_{q^{\frac{3}{2}}}^{\infty}\frac{1}{y^2}dy\nonumber\\
&\ll&\frac{\phi(q)q^{\epsilon}}{q}.
\end{eqnarray}

Similarly, we also have
\begin{eqnarray}
A_{3}=\sum_{\chi\ne \chi_{0}}\chi(a)\left(
\sum_{l=1}^{q}\frac{\bar{\chi}(l)}{l}\right)\left(
\int_{\frac{q}{a}}^{+\infty}\frac{A(\chi,y)}{y^{2}}dy
\right)=O\left(\frac{\phi(q)q^{\epsilon}}{q}\right),
\end{eqnarray}
\begin{eqnarray}
A_{4}=\sum_{\chi\ne
\chi_{0}}\chi(a)\left(\int_{\frac{q}{a}}^{+\infty}
\frac{A(\chi,y)}{y^{2}}dy\right)\left(\int_{q}^{+\infty}\frac{B(\bar{\chi},y)}{y^{2}}dy\right)=O\left(\frac{\phi(q)q^{\epsilon}}{q}\right).
\end{eqnarray}

Combining the formulas (2)-(5), we immediately obtain
$$
\sum_{\chi\ne \chi_{0}}\chi(a)\arrowvert L(1,\chi)
\arrowvert^{2m}=\frac{\phi(q)}{a}\zeta(2)\prod_{p|q}\left(1-\frac{1}{p^2}\right)+O(\log^2
q).
$$
This completes the proof of Lemma 4.

\begin{center}
\large{{\bf 3. Proof of Theorem}}
\end{center}

In this part, we will prove both of the theorems. Firstly, we will prove Theorem 1.

\textbf{Proof of Theorem 1. }According to Lemma 1 and Lemma 4, we have
\begin{eqnarray}
&&\sum_{\chi\ne\chi_0}\chi(k)\left|L(1,\chi,a)\right|^2\nonumber\\
&=&\sum_{\chi\ne\chi_0}\chi(k)\left|L(1,\chi)
-a\sum_{n=1}^{\infty}\frac{\chi(n)}{n(n+a)}\right|^2\nonumber\\
&=&\sum_{\chi\ne\chi_0}\chi(k)|L(1,\chi)|^2-a\sum_{\chi\ne\chi_{0}}\chi(k)\sum_{n=1}^{\infty}\frac{\chi(n)}{n(n+a)}L(1,\bar{\chi})-\nonumber\\
&&-a\sum_{\chi\ne
\chi_{0}}\chi(k)\sum_{n=1}^{\infty}\frac{\bar{\chi}(n)}{n(n+a)}L(1,\chi)
+a^2\sum_{\chi\ne\chi_0}\chi(k)\left|\,\sum_{n=1}^{\infty}\frac{\chi(n)}{n(n+a)}\right|^2\nonumber\\
&=&\frac{\phi(q)}{k}\zeta(2)\prod_{p|q}\left(1-\frac{1}{p^2}\right)+O(q^{\epsilon})-aM_1-aM_2+a^2M_3,\nonumber
\end{eqnarray}
where
\begin{eqnarray}
&M_1&=\sum_{\chi\ne\chi_{0}}\chi(k)\sum_{n=1}^{\infty}\frac{\chi(n)}{n(n+a)}L(1,\bar{\chi}),\nonumber\\
&M_2&=\sum_{\chi\ne
\chi_{0}}\chi(k)\sum_{n=1}^{\infty}\frac{\bar{\chi}(n)}{n(n+a)}L(1,\chi),\nonumber\\
&M_3&=\sum_{\chi\ne\chi_0}\chi(k)\left|\,\sum_{n=1}^{\infty}\frac{\chi(n)}{n(n+a)}\right|^2.\nonumber
\end{eqnarray}

Now we will estimate each term of the above.

(i) Applying Abel's identity, by analytic continuation we have
 $$
L(1,\bar{\chi})=\sum_{n=1}^{q}\frac{\bar{\chi}(n)}{n}+\int_{q}^{+\infty}\frac{B(\bar{\chi},y)}{y^2}dy,
 $$
 $$
\sum_{n=1}^{\infty}\frac{\chi(n)}{n(n+a)}=\sum_{n=1}^{N}\frac{\chi(n)}{n(n+a)}
+\int_{N}^{+\infty}\frac{(2y+a)C(\chi,y)}{y^2(y+a)^2}dy,
 $$
where $B(\chi,y)=\sum_{q\le n\le y}\chi(n)$ as defined in the proof of Lemma 4, $C(\chi,y)=\sum_{N\le n\le y}\chi(n)$, and $N>q$ is an integer.
\begin{eqnarray}
&&M_1\nonumber\\
&=&\sum_{\chi\ne
\chi_{0}}\chi(k)\sum_{n=1}^{\infty}\frac{\chi(n)}{n(n+a)}L(1,\bar{\chi})\nonumber\\
&=&\sum_{\chi\ne
\chi_{0}}\chi(k)\left(\sum_{n=1}^{N}\frac{\chi(n)}{n(n+a)}
+\int_{N}^{+\infty}\frac{(2y+a)C(\chi,y)}{y^2(y+a)^2}dy\right)
\left(\sum_{n=1}^{q}\frac{\bar{\chi}(n)}{n}+\int_{q}^{+\infty}\frac{B(\bar{\chi},y)}{y^2}dy\right)\nonumber\\
&=&\sum_{\chi\ne
\chi_{0}}\chi(k)\sum_{n=1}^{N}\frac{\chi(n)}{n(n+a)}\sum_{m=1}^{q}\frac{\bar{\chi}(m)}{m}+\sum_{\chi\ne
\chi_{0}}\chi(k)\sum_{n=1}^{N}\frac{\chi(n)}{n(n+a)}
\int_{q}^{+\infty}\frac{B(\bar{\chi},y)}{y^{2}}dy+\nonumber\\
&&+\sum_{\chi\ne
\chi_{0}}\chi(k)\sum_{n=1}^{q}\frac{\bar{\chi}(n)}{n}\int_{N}^{+\infty}\frac{(2y+a)C(\chi,y)}{y^2(y+a)^2}dy+\nonumber\\
&&+\sum_{\chi\ne
\chi_{0}}\chi(k)\int_{N}^{+\infty}\frac{(2y+a)C(\chi,y)}{y^2(y+a)^2}dy
\int_{q}^{+\infty}\frac{B(\bar{\chi},z)}{z^{2}}dz\nonumber\\
&\equiv&B_1+B_2+B_3+B_4.\nonumber
\end{eqnarray}

We will estimate each of them. Firstly we estimate $B_1$. From the
properties of Dirichlet characters and M\"{o}bius function, we
have
\begin{eqnarray}
B_1&=&\sum_{\chi\ne
\chi_{0}}\chi(k)\sum_{n=1}^{N}\frac{\chi(n)}{n(n+a)}\sum_{m=1}^{q}\frac{\bar{\chi}(m)}{m}\nonumber\\
&=&\sum_{\chi\bmod
q}\,\,\sum_{n=1}^{N}\sum_{m=1}^{q}\frac{\chi(k)\bar{\chi}(m)\chi(n)}{mn(n+a)}+O\left(\log q\right)\nonumber\\
&=&\displaystyle\mathop{\sum{'}}_{n=1}^{N}\displaystyle\mathop{\sum{'}}_{m=1}^{q}\frac{1}{mn(n+a)}\sum_{\chi\bmod
q}\chi(k)\bar{\chi}(m)\chi(n)
+O\left(\log q\right)\nonumber\\
&=&\phi(q)\displaystyle\mathop{\displaystyle\mathop{\sum{'}}_{n=1}^{N}\displaystyle\mathop{\sum{'}}_{m=1}^{q}}_{nk\equiv
m(\bmod q)}\frac{1}{mn(n+a)}
+O\left(\log q\right)\nonumber\\
&=&\phi(q)\displaystyle\mathop{\sum{'}}_{n=1}^{q}\frac{1}{n^2k(n+a)}
+\phi(q)\displaystyle\mathop{\sum{'}}_{l=1}^{N/q}\displaystyle\mathop{\sum{'}}_{n=1}^{q}\frac{1}{n(nk+lq)(n+a)}
+O\left(\log q\right)\nonumber\\
&=&\phi(q)\sum_{n=1
}^{q}\frac{1}{n^2k(n+a)}\sum_{d|(n,q)}\mu(d)+\nonumber\\
&&+O\left(\phi(q)
\displaystyle\mathop{\sum{'}}_{l=1}^{N/q}\displaystyle\mathop{\sum{'}}_{n=1}^{q}
\frac{1}{n(nk+lq)(n+a)}+\log q\right)\nonumber\\
&=&\phi(q)\sum_{d|q}\mu(d)\sum_{n=1 \atop
d|n}^{q}\frac{1}{n^2k(n+a)}+O\left(\log q\right)\nonumber\\
&=&\phi(q)\sum_{d|q}\frac{\mu(d)}{d^3}\sum_{n=1}^{q/d}\frac{1}{n^2k(n+a/d)}
+O\left(\log q\right)\nonumber\\
&=&\frac{\phi(q)}{k}\sum_{d|q}\frac{\mu(d)}{d^3}\sum_{n=1}^{\infty}\frac{1}{n^2(n+a/d)}+\nonumber\\
&&+O\left(\frac{\phi(q)}{k}\sum_{d|q}\frac{\mu(d)}{d^3}\sum_{n=q/d}^{\infty}\frac{1}{n^2(n+a/d)}
\right)
+O\left(\log q\right)\nonumber\\
&=&\frac{\phi(q)}{k}\sum_{d|q}\frac{\mu(d)}{d^3}
\left(\sum_{n=1}^{\infty}\frac{1}{\left(a/d\right)^2}\left(\frac{a/d}{n^2}+\frac{1}{n+a/d}-\frac{1}{n}\right)\right)
+O\left(\log q\right)\nonumber\\
&=&\frac{\phi(q)}{ak}\sum_{d|q}\frac{\mu(d)}{d^2}\zeta(2)
-\frac{\phi(q)}{a^2k}\sum_{d|q}\frac{\mu(d)}{d}\sum_{l=1}^{[a/d]}\frac{1}{l}
+O\left(\log q\right)\nonumber\\
&=&\frac{\phi(q)}{ak}\zeta(2)\prod_{p|q}\left(1-\frac{1}{p^2}\right)
-\frac{\phi(q)}{a^2k}\sum_{d|q}\frac{\mu(d)}{d}\sum_{l=1}^{[a/d]}\frac{1}{l}+O\left(\log
q\right).\nonumber
\end{eqnarray}

In the following we will estimate $B_2,\,B_3$ and $B_4$. According
classical estimate of Dirichlet character sums, we have
\begin{eqnarray}
|B_2|&=&\sum_{\chi\ne
\chi_{0}}\chi(k)\sum_{n=1}^{N}\frac{\chi(n)}{n(n+a)}
\int_{q}^{+\infty}\frac{B(\bar{\chi},y)}{y^{2}}dy\nonumber\\
&\le&\sum_{\chi\ne \chi_{0}}\sum_{n=1}^{N}\frac{1}{n(n+a)}
\int_{q}^{+\infty}\frac{\left|\sum_{q<n<y}\bar{\chi}(n)\right|}{y^{2}}dy\nonumber\\
&\ll&q^{\frac{1}{2}}\phi(q)\log
q \int_{q}^{\infty}\frac{1}{y^2}dy\nonumber\\
&\ll&\frac{\phi(q)\log q}{\sqrt{q}} ,\nonumber\\
|B_3|&=&\sum_{\chi\ne
\chi_{0}}\chi(k)\sum_{n=1}^{q}\frac{\bar{\chi}(n)}{n}\int_{N}^{+\infty}\frac{(2y+a)C(\chi,y)}{y^2(y+a)^2}dy\nonumber\\
&\ll&\log q\int_{N}^{\infty}\frac{2y+a}{y^2(y+a)^2}\sum_{\chi\ne
\chi_{0}}|C(\chi,y)|dy\nonumber\\
&\ll&\frac{q^\frac12\phi(q)\log^2 q}{N^2},\nonumber\\
|B_4|&=&\sum_{\chi\ne
\chi_{0}}\chi(k)\int_{N}^{+\infty}\frac{(2y+a)C(\chi,y)}{y^2(y+a)^2}dy
\int_{q}^{+\infty}\frac{B(\bar{\chi},z)}{z^{2}}dz\nonumber\\
&\le&\int_{N}^{\infty}\int_{q}^{\infty}\frac{2y+a}{y^2(y+a)^2z^2}\sum_{\chi\ne
\chi_{0}}|C(\chi,y)|\cdot|B(\bar{\chi},z)|dydz\nonumber\\
&\ll&\frac{\phi(q)\log^2 q}{N^2},\nonumber
\end{eqnarray}
where we have used the common estimates. Taking $N=q^2$, then we get
\begin{eqnarray}
M_1&=&\sum_{\chi\ne
\chi_{0}}\chi(k)\sum_{n=1}^{\infty}\frac{\chi(n)}{n(n+a)}L(1,\bar{\chi})\nonumber\\
&=&\frac{\phi(q)}{ak}\zeta(2)\prod_{p|q}\left(1-\frac{1}{p^2}\right)
-\frac{\phi(q)}{a^2k}\sum_{d|q}\frac{\mu(d)}{d}\sum_{l=1}^{[a/d]}\frac{1}{l}
+O\left(\frac{\phi(q)\log q}{\sqrt{q}}\right).\nonumber
\end{eqnarray}

This gives the asymptotic formula of $M_1$. Therefore, we have the
asymptotic formula of $M_1$.

(ii) Following the similar method in part (i), we could get the
asymptotic formula of $M_2$, that is
\begin{eqnarray}
M_2&=&\sum_{\chi\ne
\chi_{0}}\chi(k)\sum_{n=1}^{\infty}\frac{\bar{\chi}(n)}{n(n+a)}L(1,\chi)\nonumber\\
&=&\frac{\phi(q)}{ak}\zeta(2)\prod_{p|q}\left(1-\frac{1}{p^2}\right)
+\frac{\phi(q)}{a^2}\sum_{d|q}\frac{\mu(d)}{d}\sum_{l=1}^{[a/kd]}\frac{1}{l}+O\left(\log
q\right).\nonumber
\end{eqnarray}

(iii) Lastly we will derive the asymptotic formula of $M_3$. Let
$N>q$ be any integer, then we have
\begin{eqnarray}
M_3&=& \sum_{\chi\ne
\chi_{0}}\chi(k)\left|\sum_{n=1}^{\infty}\frac{\chi(n)}{n(n+a)}\right|^2\nonumber\\
&=&\sum_{\chi\ne
\chi_{0}}\chi(k)\left(\sum_{n=1}^{N}\frac{\chi(n)}{n(n+a)}
+\int_{N}^{+\infty}\frac{(2y+a)C(\chi,y)}{y^2(y+a)^2}dy\right)\times\nonumber\\
&&\times\left(\sum_{m=1}^{N}\frac{\bar{\chi}(m)}{m(m+a)}
+\int_{N}^{+\infty}\frac{(2z+a)C(\bar{\chi},z)}{z^2(z+a)^2}dz\right)\nonumber\\
&=&\sum_{\chi\ne
\chi_{0}}\chi(k)\left(\sum_{n=1}^{N}\frac{\chi(n)}{n(n+a)}
\right)\left(\sum_{m=1}^{N}\frac{\bar{\chi}(m)}{m(m+a)}
\right)+\nonumber\\
&&+\sum_{\chi\ne
\chi_{0}}\chi(k)\left(\sum_{n=1}^{N}\frac{\chi(n)}{n(n+a)}
\right)\left(\int_{N}^{+\infty}\frac{(2z+a)C(\bar{\chi},z)}{z^2(z+a)^2}dz\right)+\nonumber\\
&&+\sum_{\chi\ne
\chi_{0}}\chi(k)\left(\int_{N}^{+\infty}\frac{(2y+a)C(\chi,y)}{y^2(y+a)^2}dy\right)
\left(\sum_{m=1}^{N}\frac{\bar{\chi}(m)}{m(m+a)}\right)+\nonumber\\
&&+\sum_{\chi\ne
\chi_{0}}\chi(k)\left(\int_{N}^{+\infty}\frac{(2y+a)C(\chi,y)}{y^2(y+a)^2}dy\right)
\left(\int_{N}^{+\infty}\frac{(2z+a)C(\bar{\chi},z)}{z^2(z+a)^2}dz\right)\nonumber\\
&=&\sum_{n=1}^{N}\frac{1}{n(n+a)}
\sum_{m=1}^{N}\frac{1}{m(m+a)}\left(\sum_{\chi\ne
\chi_{0}}\chi(nk)\bar{\chi}(m)
\right)+\nonumber\\
&&+O\left(\int_{N}^{+\infty}\frac{(2y+a)\sum_{\chi\ne \chi_0}|C(\chi,y)|}{y^2(y+a)^2}dy\right)\nonumber\\
&=&\phi(q)\displaystyle\mathop{\displaystyle\mathop{\sum{'}}_{n=1}^{N}\displaystyle\mathop{\sum{'}}_{m=1}^{N}}_{m\equiv
nk(\bmod
q)}\frac{1}{mn(m+a)(n+a)}+O(1)+O\left(\frac{q^{\frac12}\phi(q)\log q}{N^2}\right)\nonumber\\
&=&\phi(q)\displaystyle\mathop{\sum{'}}_{n=1}^{N}\frac{1}{n^2k(n+a)(nk+a)}
+2\phi(q)\displaystyle\mathop{\displaystyle\mathop{\sum{'}}_{n=1}^{N}\displaystyle\mathop{\sum{'}}_{m=1}^{N}}_{m\equiv
nk(\bmod q),\,m>nk}\frac{1}{mn(m+a)(n+a)}+O(1)\nonumber\\
&=&\frac{\phi(q)}{k^2}\displaystyle\mathop{\sum{'}}_{n=1}^{N}\frac{1}{n^2(n+a)(n+a/k)}\nonumber\\
&&+O\left(\phi(q)\displaystyle\mathop{\sum{'}}_{l=1}^{N/q}
\displaystyle\mathop{\sum{'}}_{n=1}^{N}\frac{1}{n(n+a)(lq+nk)(lq+nk+a)}\right)
+O(1)\nonumber\\
&=&\frac{\phi(q)}{k^2}\sum_{n=1}^{N}\frac{1}{n^2(n+a)(n+a/k)}\sum_{d|(n,q)}\mu(d)+O\left(\frac{\phi(q)\log
N}{q}\right)
+O(1)\nonumber\\
&=&\frac{\phi(q)}{k^2}\sum_{d|q}\frac{\mu(d)}{d^4}\sum_{n=1}^{N/d}\frac{1}{n^2(n+a/d)(n+a/kd)}
+O\left(1\right)\nonumber\\
&=&\frac{\phi(q)}{k^2}\sum_{d|q}\frac{\mu(d)}{d^4}\sum_{n=1}^{\infty}\frac{1}{n^2(n+a/d)(n+a/kd)}+\nonumber\\
&&+O\left(\frac{\phi(q)}{k^2}\sum_{d|q}\frac{\mu(d)}{d^4}\sum_{n=N/d}^{\infty}\frac{1}{n^2(n+a/d)(n+a/kd)}\right)
+O\left(1\right)\nonumber\\
&=&\frac{\phi(q)}{k^2}\sum_{d|q}\frac{\mu(d)}{d^4}
\sum_{n=1}^{\infty}\frac{d^3}{a^3}\left(\frac{ka/d}{n^2}-\frac{k(k+1)}{n}-\frac{k/(k-1)}{n+a/d}+\frac{k^3/(k-1)}{n+a/kd}\right)\nonumber\\
&&+O\left(1\right)\nonumber\\
&=&\frac{\phi(q)}{a^2k}\zeta(2)\prod_{p|q}\left(1-\frac{1}{p^2}\right)+
\frac{\phi(q)}{a^3k(k-1)}\sum_{d|q}\frac{\mu(d)}{d}\sum_{l=1}^{[\frac{a}{d}]}\frac{1}{l}+\nonumber\\
&&
-\frac{k\phi(q)}{a^3(k-1)}\sum_{d|q}\frac{\mu(d)}{d}\sum_{l=1}^{[\frac{a}{kd}]}\frac{1}{l}+O\left(1\right).\nonumber
\end{eqnarray}
Taking $N=q^2$, we have
\begin{eqnarray}
M_3&=&\sum_{\chi\ne
\chi_{0}}\chi(k)\left|\sum_{n=1}^{\infty}\frac{\chi(n)}{n(n+a)}\right|^2\nonumber\\
&=&\frac{\phi(q)}{a^2k}\zeta(2)\prod_{p|q}\left(1-\frac{1}{p^2}\right)+\frac{\phi(q)}{a^3k(k-1)}\sum_{d|q}\frac{\mu(d)}{d}\sum_{l=1}^{[\frac{a}{d}]}\frac{1}{l}-\nonumber\\
&&
-\frac{k\phi(q)}{a^3(k-1)}\sum_{d|q}\frac{\mu(d)}{d}\sum_{l=1}^{[\frac{a}{kd}]}\frac{1}{l}+O\left(1\right).\nonumber
\end{eqnarray}

Combining the estimates of (i),(ii), (iii) and Lemma 4, we
immediately obtain
\begin{eqnarray}
&&\sum_{\chi\ne\chi_0}\chi(k)\left|L(1,\chi,a)\right|^2\nonumber\\
&=&\sum_{\chi\ne\chi_0}\chi(k)\left|L(1,\chi)
-a\sum_{n=1}^{\infty}\frac{\chi(n)}{n(n+a)}\right|^2\nonumber\\
&=&\sum_{\chi\ne\chi_0}\chi(k)|L(1,\chi)|^2-a\sum_{\chi\ne
\chi_{0}}\chi(k)\sum_{n=1}^{\infty}\frac{\bar{\chi}(n)}{n(n+a)}L(1,\chi)-\nonumber\\
&&-a\sum_{\chi\ne\chi_{0}}\chi(k)\sum_{n=1}^{\infty}\frac{\chi(n)}{n(n+a)}L(1,\bar{\chi})
+a^2\sum_{\chi\ne\chi_0}\left|\,\sum_{n=1}^{\infty}\frac{\chi(n)}{n(n+a)}\right|^2\nonumber\\
&=&\frac{\phi(q)}{a(k-1)}\sum_{d|q}\frac{\mu(d)}{d}\sum_{l=1}^{[\frac{a}{d}]}\frac{1}{l}
-\frac{\phi(q)}{a(k-1)}\sum_{d|q}\frac{\mu(d)}{d}\sum_{l=1}^{[\frac{a}{kd}]}\frac{1}{l}
+O\left(\frac{\phi(q)\log q}{\sqrt{q}}\right).\nonumber\\
&=&\frac{\phi(q)}{a(k-1)}\sum_{d|q}\frac{\mu(d)}{d}\sum_{l=[\frac{a}{kd}]+1}^{[\frac{a}{d}]}\frac{1}{l}
+O\left(\frac{\phi(q)\log q}{\sqrt{q}}\right),\nonumber
\end{eqnarray}
where the $O$ constant depends on $a,k$. This proves Theorem 1.\\

Next, we shall complete the proof of Theorem 2. Theorem 1 will be useful in the proof.

\textbf{Proof of Theorem 2.}
Firstly from Lemma 2, for any $p\ge 3$, we have
\begin{eqnarray}
&&\sum_{\chi\ne\chi_0}\left|\sum_{x=1}^{p-1}\chi(x)e\left(\frac{f(x)}{p}\right)\right|^2
|L(1,\chi,a)|^{2}\nonumber\\
&=&\sum_{\chi\ne\chi_0}\left(p-1+
\sum_{x=2}^{p-1}\chi(x)\sum_{y=1}^{p-1}e\left(\frac{g(y,x)}{p}\right)\right)
|L(1,\chi,a)|^{2}\nonumber\\
&=&(p-1)\sum_{\chi\ne\chi_0}|L(1,\chi,a)|^2+
\sum_{x=2}^{p-1}\sum_{y=1}^{p-1}e\left(\frac{g(y,x)}{p}\right)\sum_{\chi\ne\chi_0}\chi(x)|L(1,\chi,a)|^{2}\nonumber\\
&=&(p-1)\sum_{\chi\ne\chi_0}\arrowvert
L(1,\chi,a)\arrowvert^{2}\nonumber\\
&&
+\sum_{x=2}^{p-1^{*}}\sum_{y=1}^{p-1}e\left(\frac{g(y,x)}{p}\right)\sum_{\chi\ne\chi_0}\chi(x)
\arrowvert
L(1,\chi,a)\arrowvert^{2}\nonumber\\
&&+\sum_{x=2}^{p-1^{**}}\sum_{y=1}^{p-1}e\left(\frac{g(y,x)}{p}\right)\sum_{\chi\ne\chi_0}\chi(x)
\arrowvert L(1,\chi,a)\arrowvert^{2},\nonumber
\end{eqnarray}
where $g(y,x)=\sum_{i=0}^{k}a_{i}(x^{i}-1)y^{i}$, $b_i=a_i(x^i-1)$
and $\sum_{x=2}^{p-1^{*}}\sum_{y=1}^{p-1}$ and
$\sum_{x=2}^{p-1^{**}}\sum_{y=1}^{p-1}$ means
$p\dag(b_0,b_1,\dots,b_k)$ and $p|(b_0,b_1,\dots,b_k)$
respectively. Then we will estimate the two sums respectively.

  (1) When $p\dag(b_0,b_1,\dots,b_k)$, according to Lemma 3 and Theorem 1, we have \\
\begin{eqnarray}
&&\left|\sum_{x=2}^{p-1^{*}}\sum_{y=1}^{p-1}e\left(\frac{g(y,x)}{p}\right)\sum_{\chi\ne\chi_0}\chi(x)
\arrowvert
L(1,\chi,a)\arrowvert^{2}\right|\nonumber\\
&\ll&\sum_{x=2}^{p-1}p^{1-\frac{1}{k}+\epsilon}\left|\sum_{\chi\ne\chi_0}\chi(x)
\arrowvert L(1,\chi,a)\arrowvert^{2}\right|\nonumber\\
&\ll&\sum_{x=2}^{p-1}\frac{p^{2-\frac{1}{k}+\epsilon}\log p}{a(x-1)}\nonumber\\
&\ll&p^{2-\frac{1}{k}+\epsilon}.\nonumber
\end{eqnarray}

  (2) When $p\mid (b_0,b_1,\dots,b_k),$ \,i.e.\,$p\mid b_0,p\mid b_1,\dots,p\mid b_k,$
  since $p\dag(a_0,a_1,\dots,a_k),$ there is at
least one $a_{l}$ such that $p\dag a_{l}$, then for this $l$ we
must have $p\mid (x^{l}-1)$,i.e. $x^l\equiv 1(\mbox{mod}p)$. But
in the set $\{2,3,\dots,p-1\}$, there are at most $l-1$ numbers
$x$ such that $p\mid (x^l-1)$. Also $l-1<l\le k$, $x^l>x^l-1\ge
p$, so $x>p^{\frac{1}{l}}\ge p^{\frac{1}{k}}$.
Then from Lemma 3 and Theorem 1 we have\\
\begin{eqnarray}
&&\left|\sum_{x=2}^{p-1^{**}}\sum_{y=1}^{p-1}e\left(\frac{g(y,x)}{p}\right)\sum_{\chi\ne\chi_0}\chi(x)
\arrowvert L(1,\chi,a)\arrowvert^{2}\right|\nonumber\\
&\le&\sum_{x=2}^{p-1^{**}}\left|
\sum_{y=1}^{p-1}e\left(\frac{g(y,x)}{p}\right)\right|\left|\sum_{\chi\ne\chi_0}\chi(x)
\arrowvert L(1,\chi,a)\arrowvert^{2}\right|\nonumber\\
&\ll&\sum_{x=2}^{p-1^{**}}(p-1)\frac{p\log p}{a(x-1)}\nonumber\\
&\ll&p^2\log p\left(\max_{p^{1/k}\le x<
p}\left(\frac{1}{x-1}\right)\right)
\times\sharp\{x:x\in\{2,3,\dots,p-1\}\quad\mbox{with}\quad x^{l}\equiv 1(\bmod\,\, p)\}\nonumber\\
&\ll&kp^{2-\frac{1}{k}+\epsilon}.\nonumber
\end{eqnarray}

Therefore,\,\,combining (1),\,(2), the formula (1) and Theorem 1, we
get the asymptotic formula
\begin{eqnarray}
&&\sum_{\chi\ne\chi_0}\left|\sum_{x=1}^{p-1}\chi(x)e\left(\frac{f(x)}{p}\right)\right|^2
|L(1,\chi,a)|^{2}\nonumber\\
&=&(p-1)\sum_{\chi\ne\chi_0}\arrowvert
L(1,\chi,a)\arrowvert^{2} \nonumber\\
&&+O\left(\left|\sum_{x=2}^{p-1^{*}}\sum_{y=1}^{p-1}e\left(\frac{g(y,x)}{p}\right)\sum_{\chi\ne\chi_0}\chi(x)
\arrowvert
L(1,\chi,a)\arrowvert^{2}\right|\right)\nonumber\\
&&+O\left(\left|\sum_{x=2}^{p-1^{**}}\sum_{y=1}^{p-1}e\left(\frac{g(y,x)}{p}\right)\sum_{\chi\ne\chi_0}\chi(x)
\arrowvert
L(1,\chi,a)\arrowvert^{2}\right|\right)\nonumber\\
&=&p^2\sum_{d|q}\frac{\mu(d)}{d^2}\zeta\left(2,\frac{a}{d}\right)
-\frac{4p^2}{a}\sum_{d|q}\frac{\mu(d)}{d}\sum_{l=1}^{[\frac{a}{d}]}\frac{1}{l}
+O(p^{2-\frac{1}{k}+\epsilon}),\nonumber
\end{eqnarray}
where $\mu$ is the M\"{o}bius function and
$\zeta(s,\alpha)(s=\sigma+it,\,\alpha>0)$ is the Hurwitz zeta
function, the $O$ constant is depending on $k$,$a$ and $\epsilon$.
This completes the proof of Theorem 2.\\


\begin{thebibliography}{99}
\bibitem{S.H}
D. R. Heath-Brown, An asymptotic series for the mean value of
Dirichlet $L$-functions, Comment. Math. Helvetici 56, 1981,
148-161.
\bibitem{Z.Y}
W. P. Zhang, On the second mean value of Dirichlet $L$-functions,
Chinese Annals of Mathematics A 11, 1990, 121-127.(in Chinese)
\bibitem{S.H}
R. Balasubramanian, A note on Dirichlet $L$-functions, Acta Arith.
38, 1980, 273-283.
\bibitem{Z.Y}
W. P. Zhang, Y. Yi and X. L. He, On the $2k$-th Power Mean of
Dirichlet $L$-Functions With the Weight of General Kloosterman
Sums, Journal of Number Theory, 84, 2000, 199-213.
\bibitem{S.H}
Y. Yi and W. P. Zhang, On the $2k$-th Power Mean of Dirichlet
$L$-Functions With the Weight of Gauss Sums, Advances in
Mathematics 31(6), 2002, 517-526.
\bibitem{Y.Z}
R. Ma, J. H. Zhang, Y. L. Zhang, On the $2m$-th power mean of
Dirichlet $L$-functions with the weight of trigonometric sums, Proc.
Indian Acad. Sci. (Math. Sci.) 119(4), 2009, 411-421.
\bibitem{L.K}
B. C. Berndt, Generalized Dirichlet series and Hecke's functional
equation, Proc.Edinburgh Math. Soc. 15(2), 1967, 309-313.
\bibitem{S.H}
B. C. Berndt, Identities involving the coefficients of a class of
Dirichlet series. III, Trans.Amer. Math. Soc. 146, 1969, 323-342.
\bibitem{S.H}
B. C. Berndt, Identities involving the coefficients of a class of
Dirichlet series. IV, Trans. Amer. Math. Soc. 149, 1970, 179-185.
\bibitem{S.H}
R. Ma, Y. Yi, Y. L. Zhang, On the mean value of the gerneralized
Dirichlet $L$-functions, Czechoslovak Mathematical Journal 60(135)
, 2010, 597-620.
\bibitem{S.H}
L. J. Mordell, On a sum analogous to a Gauss's sum, Quart. J.
Math., Oxford 3, 1932, 161-167.
\bibitem{L.C}
L. Carlitz and S. Uchiyama, Bounds for exponential sums, Duke
Math. J. 24(1), 1957, 37-41.
\bibitem{A.W}
A. Weil, On some exponential sums, Proc. nat. Acad, Sci, USA 34, 1948,
204-207.
\bibitem{H.D}
H. Davenport, On certain exponential sums, J. Reine Angew.
Math. 169, 1933, 158-176.
\bibitem{L.K}
L. K. Hua, On exponential sums over an algebraic field, Canadian
J. Math. 3, 1951, 44-51.
\bibitem{T.M}
M. A. Tom, Introduction to Analytic Number Theory,
Springer-Verlag, New York, 1976.
\bibitem{C.D}
C. D. Pan and C. B. Pan, Element of the Analytic Number Theory,
Science Press, Beijing, 1991 (in Chinese).
\end{thebibliography}
\end{document}